\documentclass[twoside]{article}

 \usepackage[dvips]{graphicx,epsfig}

\input amssym.def
\input amssym.tex
\usepackage{bbm}

\skip\footins 22pt plus 2pt minus 2pt
\newdimen\jot
\def\mqth{\mathsurround=0pt }
\dimendef\dimenq=0
\def\openup{\afterassignment\qpenup\dimenq=}
\def\qpenup{\advance\lineskip\dimenq
  \advance\baselineskip\dimenq \advance\lineskiplimit\dimenq}
\def\eqalign#1{\,\vcenter{\openup1\jot \mqth
  \ialign{\strut\hfil$\displaystyle{##}$&$\displaystyle{{}##}$\hfil
  \crcr#1\crcr}}\,}
\newif\ifdtqp
\def\displqy{\global\dtqptrue \openup1\jot \mqth
  \everycr{\noalign{\ifdtqp \global\dtqpfalse
     \vskip-\lineskiplimit \vskip\normallineskiplimit
     \else \penalty\interdisplaylinepenalty \fi}}}
\def\displaylines#1{\displqy
  \halign{\hbox to\displaywidth{$\hfil\displaystyle##\hfil$}\crcr
  #1\crcr}}
\newskip\centerinq \centerinq=0pt plus 1000pt minus 1000pt

\def\bah#1{\overline#1}

\def\dast{{\displaystyle{\ast}}}

\newcommand{\ph}{\varphi}
\newcommand{\reals}{\mathbbm{R}}

\newcommand{\Harm}{\mathop{\rm Harm}\nolimits}
\newcommand{\Charm}{\mathop{\rm CHarm}\nolimits}
\newcommand{\Eharm}{\mathop{\rm EHarm}\nolimits}
\renewcommand{\tan}{{\mathop{\rm tan}\nolimits}}
\newcommand{\norm}{{\mathop{\rm norm}\nolimits}}
\newcommand{\im}{\mathop{\rm im}\nolimits}
\renewcommand{\ker}{\mathop{\rm ker}\nolimits}

\begin{document}
\thispagestyle{empty}
\large
\addtolength{\baselineskip}{2pt}
\centerline{\Large{\bf Cohomology of Harmonic Forms on Riemannian 
Manifolds With Boundary\rm}}
\medskip
\centerline{Sylvain Cappell, Dennis DeTurck, Herman Gluck, 
and Edward Y.\ Miller}
\medskip
\centerline{\it To Julius Shaneson on the occasion of his \rm 60\it th
birthday\rm}

\medskip
\noindent\bf 1. Introduction\rm

The main result of this article is the following.

\medskip
\noindent\bf Theorem 1\it.\ \ Let $M$ be a compact, 
connected, oriented, smooth Riemannian n-dimensional manifold with
non-empty
boundary. Then the cohomology of the complex $(\Harm^\dast(M),d)$ 
of harmonic forms on $M$ is given by the direct sum:
$$
H^p(\Harm^\dast(M),d) \cong H^p(M;\reals) + H^{p-1}(M;\reals)
$$
for $p=0,1,\ldots, n$\rm.

\medskip
Let $M$ be a smooth compact $n$-manifold, and $\Omega^\dast(M)$ the space
of smooth differential forms on $M$. The classical theorem of 
de\,Rham [1931] asserts that the cohomology of the complex
$$0\longrightarrow \Omega^0(M) 
{\smash{\mathop{\longrightarrow}\limits^{d}}}  \Omega^1(M)
{\smash{\mathop{\longrightarrow}\limits^{d}}}  \Omega^2(M)
{\smash{\mathop{\longrightarrow}\limits^{d}}}  \cdots
{\smash{\mathop{\longrightarrow}\limits^{d}}}  \Omega^{n-1}(M)
{\smash{\mathop{\longrightarrow}\limits^{d}}}  \Omega^n(M)
\longrightarrow 0,$$
where $d$ is exterior differentiation,
is isomorphic to the cohomology of $M$ with real coefficients.
In other words,
$${\ker d:\Omega^p(M)\to \Omega^{p+1}(M)\over
\im d:\Omega^{p-1}(M)\to \Omega^p(M)} \cong H^p(M;\reals).$$

If $M$ is oriented,
a Riemannian metric on $M$ gives rise to an $L^2$ inner product
$$\langle \alpha,\beta\rangle = \int_M \alpha\wedge\dast\beta$$ 
on $\Omega^\dast(M)$, where $\dast$ denotes
the Hodge star operator, and to the co-differential 
$$\delta=(-1)^{n(p+1)+1}\dast d\dast\colon\Omega^p(M)\to
\Omega^{p-1}(M),$$ 
which on a closed manifold
is the $L^2$-adjoint of the exterior differential $d$. As usual,
one defines the Laplacian by 
$$\Delta=d\delta + \delta d\colon\Omega^p(M)\to\Omega^p(M),$$ 
and  \it harmonic \rm differential forms $\omega$
as those satisfying $\Delta\omega=0$. 


The exterior differential $d$, since it commutes with $\Delta$,
preserves harmonicity of forms, and hence $(\Harm^\dast(M),d)$ 
is a subcomplex 
of the de\,Rham complex $(\Omega^\dast(M),d)$. It is therefore
natural to compute the cohomology of 
this complex, which we call the \it harmonic cohomology 
\rm of $M$.

When $M$ is a closed manifold, a form $\omega$ is harmonic if and only
if it is both closed ($d \omega =0$) and co-closed
($\delta \omega =0$). 
In this case, all of the maps in the complex $(\Harm^\dast(M),d)$
are zero, and so
$$
H^p(\Harm^\dast(M),d) = \Harm^p(M) \cong H^p(M;\reals),
$$
according to the classical theorem of Hodge [1933].

By contrast, when $M$ is connected and 
has non-empty boundary, it is possible
for a\\ 
$p$-form to be harmonic without being both
closed and co-closed. Some of these, which are exact (that is,
in the image of $d$), although not the
exterior derivatives of \it harmonic \rm $p-1$-forms, represent
the ``echo'' of the ordinary $p-1$-dimensional cohomology within
the $p$-dimensional harmonic cohomology reported in Theorem 1.


\noindent\bf Example\rm: 
%
\ Let $M$ be the annulus $a^2\le x^2+y^2\le b^2$ in the
$xy$-plane
$\reals^2$. Then the 2-form $\omega=-{1\over 2}\log(x^2+y^2)\,dx\wedge dy$
is
harmonic because $\log(x^2+y^2)$ is a harmonic function on the annulus. 
Even though $\omega$ is exact, it is not the exterior derivative of any
harmonic 
1-form. Hence it represents a nonzero element of the 2-dimensional
harmonic cohomology $H^2(\Harm^\dast(M),d)$. 
Also, $\omega$ is not co-closed, in fact 
$$\delta\omega=\ph=
(-y\,dx+x\,dy)/(x^2+y^2),$$ which represents a generator of the
1-dimensional
cohomology $H^1(M;\reals)\cong \reals$. The equation $\delta\omega=\ph$ is
the signal that the 2-form $\omega$ is the echo of the 1-form $\ph$, as
we will see in Lemma 3 below. 

\medskip

\noindent\bf Remarks\rm: 1.\ The co-differential $\delta$ also commutes with
the Laplacian $\Delta$, and therefore $(\Harm^\dast(M),\delta)$ is a 
subcomplex of $(\Omega^\dast(M),\delta)$. We can apply the Hodge
star operator to the isomorphism given by Theorem 1 and
replace $n-p$ by $p$ to obtain
$$
H^p(\Harm^\dast(M),\delta) \cong 
H^p(M, \partial M;\reals) +  H^{p+1}(M,\partial M;\reals),
$$
where $H^p(M,\partial M;\reals)$ is the cohomology 
of $M$ relative to its boundary. 
In this case, the homological ``echo'' is  shifted 
\it down \rm by one unit.

\noindent 2.\  Theorem 1, Remark 1 and
their proofs can be readily
generalized to harmonic forms with coefficients in a flat bundle
with metric. 

\noindent 3.\ It would be interesting to understand to what 
extent these results have analogues
for harmonic forms on smoothly stratified manifolds
with singularities. In this regard, we
look to the work of Cheeger [1980], and Mazzeo and Melrose [1999].
It would also
be interesting to consider analogous questions for the
$\overline{\partial}$
operator, cf.\ Epstein [2005] for some results on solutions of
$\bah{\partial}$ equations with modified $\bah{\partial}$-Neumann
conditions along the boundary. 

\noindent 4.\ In the proof of Theorem 1 given here, we never
make use of the structure of the space of all harmonic forms on $M$,
but focus only on the closed ones. Theorem 3.25 on pp 48--49 of
Parsley [2004] gives the structure of all harmonic vector fields
on a compact Riemannian 3-manifold with boundary; the analogous
result holds for differential $p$-forms on compact Riemannian $n$-manifolds
with boundary.

\bigskip

\addtolength{\baselineskip}{0pt}
  
\noindent \bf 2. The Hodge Decomposition Theorem\rm

\noindent\bf Conventions and definitions.\rm \ \ In what
follows, the reference to the manifold $M$ is understood,
and so we omit it and write $\Omega^p$ for the space of smooth
differential
$p$-forms on $M$. We will write $C^p$ and $cC^p$ for the spaces of closed
and co-closed
$p$-forms on $M$, and $E^p$ and $cE^p$ for the spaces of exact and
co-exact (that is, in the image of $\delta$) $p$-forms on $M$. 
We juxtapose letters to indicate intersections
of spaces, so $CcC^p$ is the subspace of $p$-forms which are both closed
and co-closed (these were called \it harmonic fields \rm by Kodaira
[1949]). Similarly, $EcC^p=E^p\cap cC^p\subset CcC^p$ and 
$CcE^p=C^p\cap cE^p\subset CcC^p$.
Finally, we use the symbol $+$ between spaces to indicate a direct
sum, and reserve $\oplus$ for an orthogonal direct sum.

To prepare for the proof of Theorem 1, we consider boundary conditions
on differential forms and the related Hodge decompositions of 
$\Omega^\dast$ on manifolds with boundary.
Along the boundary of $M$,
any smooth differential $p$-form $\omega$ has a natural 
decomposition into tangential and normal components. For
$x\in \partial M$, we write 
$$
\omega(x)= \omega_{\tan}(x) + \omega_{\norm}(x),
$$
where $\omega_{\tan}(x)$ agrees with 
$\omega(x)$ when evaluated on a $p$-tuple of vectors, all of 
which are tangent to $\partial M$, but is zero if any one
of the vectors is orthogonal to $\partial M$. We then define 
$\omega_{\norm}(x)$ by the above equation. We have that
$\omega_{\tan}(x)=0$ if and only if the restriction
$(\omega\vert_{\partial M})(x)=0$.

Let $\Omega^p_N$ be the space of smooth $p$-forms on $M$ that satisfy
\it Neumann boundary conditions \rm at every point of $\partial M$,
$$\Omega^p_N=\{\omega\in \Omega^p\, |\, \omega_\norm
=0\},$$
and similarly let $\Omega^p_D$ be the space of smooth 
$p$-forms on $M$ that satisfy \it Dirichlet boundary conditions \rm
at every point of $\partial M$,
$$\Omega^p_D=\{\omega\in \Omega^p\, |\, \omega_\tan
=0\}.$$
We define
$cE^p_N=\delta(\Omega^{p+1}_N)$ and $E^p_D=d(\Omega^{p-1}_D)$,
and emphasize that the boundary conditions are applied \it before \rm
we take co-differentials and differentials.

As noted above, on a closed manifold, 
$CcC^p(M)$ and
$\Harm^p(M)$ coincide, but in the presence of a boundary,
there are more harmonic forms than fields. We apply the boundary
conditions to $CcC^p(M)$ as follows:
$$\eqalign{
CcC^p_N&=\{\omega\in \Omega^p\,|\,d
\omega=0,\ \delta\omega=0,\ \omega_\norm=0\}\cr
CcC^p_D&=\{\omega\in \Omega^p\,|\,d
\omega=0,\ \delta\omega=0,\ \omega_\tan=0\}.
}$$

\medskip

\noindent\bf Hodge Decomposition Theorem\footnote{The Hodge Decomposition 
Theorem arose historically with increasing
generality in the
papers and books of de\,Rham [1931], Hodge [1933], Weyl [1940],
Hodge [1941], Tucker
[1941],
Weyl [1943], Bidal and de\,Rham [1946], Kodaira [1949],
Duff [1952], Duff and Spencer [1952], de\,Rham [1955], Friedrichs [1955],
Conner [1955]
and Morrey [1956].
}
\rm.\ \it Let $M$ be 
a compact, connected, oriented, smooth Riemannian $n$-manifold, 
with or without boundary. Then we have the orthogonal direct
sum
$$\Omega^p=cE^p_N\oplus CcC^p \oplus E^p_D.\eqno{(1)}$$
Furthermore, 
$$CcC^p=CcC^p_N\oplus EcC^p = CcE^p \oplus CcC^p_D.$$
\rm

\medskip

When the manifold $M$ is closed, the boundary conditions are vacuous,
and we get the original Hodge decomposition, 
\ \ $\Omega^p=cE^p\oplus CcC^p\oplus E^p$.
In this case, $C^p=CcC^p\oplus E^p$, and thus
$CcC^p$ is
the orthogonal complement of the exact 
$p$-forms within the closed
ones, so $CcC^p\cong H^p(M;\reals)$. Likewise, $cC^p=cE^p\oplus CcC^p$,
and so $CcC^p$ is simultaneously the
orthogonal complement of the co-exact $p$-forms within the co-closed
ones. 

When the boundary of $M$ is non-empty, the space $C^p$ of
closed $p$-forms decomposes as
$$C^p = CcC^p\oplus E^p_D =CcC^p_N\oplus EcC^p
\oplus E^p_D =CcC^p_N\oplus E^p.$$
Thus, $CcC^p_N$ is the orthogonal complement of the exact $p$-forms 
within the closed ones, so $CcC^p_N\cong H^p(M;\reals)$. Similarly,
the space $cC^p$ of co-closed $p$-forms decomposes as
$$cC^p = cE^p_N \oplus CcC^p =
cE^p_N\oplus CcE^p \oplus CcC^p_D
 =cE^p\oplus CcC^p_D.$$
Thus, $CcC^p_D$ is the orthogonal complement of the co-exact $p$-forms
within the co-closed ones, so $CcC^p_D\cong H^p(M,\partial M;\reals)$.

All the decompositions given above
are canonical, once the Riemannian metric on $M$ is specified.

\vfill
\eject

\noindent\bf 3. The image of the Laplacian \rm

If $M$ is a closed, oriented Riemannian $n$-manifold, the Hodge Decomposition
Theorem tells us that
$\Omega^p = cE^p \oplus CcC^p \oplus E^p$.
The Laplacian $\Delta$ acting on $p$-forms is self-adjoint, and its image
$\Delta(\Omega^p)$ is the orthogonal complement $cE^p \oplus E^p$
of its kernel $CcC^p$. Thus
$\Omega^p = CcC^p\oplus \Delta(\Omega^p)$. 

By contrast, when the boundary of the manifold is non-empty, we have

\medskip

\noindent\bf Lemma 1\it.\ \ Let $M$ be a compact, connected, oriented,
smooth Riemannian $n$-manifold with non-empty boundary.
Then the Laplacian on forms, 
$\Delta\colon \Omega^p\to \Omega^p$,
is surjective.\rm

\medskip
\noindent\it Proof\rm.\ \ Equation (1) in the Hodge Decomposition Theorem
asserts that\\
$\Omega^p= cE^p_N \oplus CcC^p \oplus E^p_D$,
and we will compute the image of the Laplacian on each summand.

On $cE^p_N$, we have $\Delta=\delta d$. Since
$C^p=CcC^p\oplus E^p_D$,
the exterior derivative $d$ must take $cE^p_N$ isomorphically to
$E^{p+1}=EcC^{p+1}\oplus E^{p+1}_D$. Applying the
co-differential
$\delta$ to this, we see that $\delta$ kills $EcC^{p+1}$
and takes
$E^{p+1}_D$ isomorphically to $cE^p$. Thus
$$\Delta(cE^p_N) = cE^p = cE^p_N\oplus CcE^p.$$
Likewise,
$$\Delta(E^p_D) = E^p = EcC^p\oplus E^p_D.$$
And naturally, $\Delta(CcC^p)=0$.

Referring again to the Hodge decomposition (1), we see that the only way
that the Laplacian $\Delta\colon \Omega^p\to \Omega^p$ could fail to be
surjective would be
for  $CcE^p$ and $EcC^p$ to fail
to span $CcC^p$. But from the Hodge Decomposition Theorem, the orthogonal 
complement of $CcE^p$ in $CcC^p$ is $CcC^p_D\cong H^p(M,\partial
M;\reals)$, and the orthogonal complement of $EcC^p$ in $CcC^p$
is  $CcC^p_N\cong H^p(M;\reals)$. Thus the subspaces in question both have 
finite codimension in $CcC^p$, and so the only way they could fail to span
$CcC^p$ 
would be for some non-zero $\omega\in CcC^p$ to be orthogonal
to both subspaces. This would force $\omega$ to lie in $CcC^p_D\cap
CcC^p_N$,
telling
us that $\omega$ is closed, co-closed, and vanishes on the boundary of $M$. 
But such a form must be zero, according to the following Lemma, which will complete
the proof of Lemma 1.

\medskip

\noindent\bf Lemma 2\it.\ \ On a connected, oriented, smooth Riemannian 
$n$-manifold with non-empty boundary, a smooth differential form which is both
closed and co-closed, and which vanishes on the boundary, must be 
identically zero.\rm

\medskip

\noindent In order to prove Lemma 2, we will appeal to the ``strong 
unique continuation theorem'', orginally
due to Aronszajn [1957], Aronszajn, Krzywicki and Szarski [1962],
and given by Kazdan [1988]
in the following form:
  
\medskip

\noindent\bf Strong Unique Continutation Theorem\rm.
\it
Let $N$ be a Riemannian manifold with Lipshitz continuous
metric, and let $\omega$ be a differential 
form having first derivatives
in $L^2$ that satisfies $\Delta \omega =0$. If $\omega$ has a zero of
infinite order at some point in $N$, then $\omega$ is
identically zero. \rm

\medskip

\noindent\it Proof of Lemma \rm 2.\ \ Let $M$ be a connected, oriented,
smooth Riemannian $n$-manifold with non-empty boundary,
and $\omega$ a smooth differential $p$-form on $M$ which 
is closed, co-closed, and vanishes on $\partial M$. We will show
that $\omega$ is identically zero. Since the
result is local, we can take $M$ to be the upper half-space
in $\reals^n$, with $\partial M=\reals^{n-1}$.

Extend the metric from the upper half-space to all of $\reals^n$
by reflection in $\reals^{n-1}$. The resulting metric will be 
Lipschitz continuous. Extend the $p$-form $\omega$ to all of 
$\reals^n$ by making it odd with respect to reflection in 
$\reals^{n-1}$. Because the original $\omega$
vanished on $\reals^{n-1}$ and was closed and co-closed, the extended
$\omega$ will be of class $C^1$ and will be closed and co-closed on
all of $\reals^n$.

These facts, together
with the vanishing of $\omega$ on $\reals^{n-1}$, are enough to 
show that the first derivatives of the coefficients of
$\omega$ vanish along $\reals^{n-1}$, even when
computed in the normal direction.
Repeated differentiation of the equations which express the fact that
$\omega$ is closed and co-closed, together with the vanishing of
$\omega$ on $\reals^{n-1}$, show that all higher partial derivatives of
the coefficients of $\omega$ vanish on $\reals^{n-1}$.
In 
other words, $\omega$ vanishes to infinite order at each
point of $\reals^{n-1}$. 

The Strong Unique Continuation Theorem then implies that 
$\omega$ must be identically zero on all of $\reals^n$. 
Since $M$ was assumed to be connected, $\omega$ must 
be identically zero on all of $M$. This completes the proof
of Lemma 2, and with it, the proof of Lemma 1.

For a different proof of Lemma 1, see Theorem 3.4.10 on page
137 of Schwarz [1995]. 

\bigskip

\noindent\bf 4. Proof of Theorem 1\rm

To prove Theorem 1, we must show that
$$H^p(\Harm^\dast(M),d)\cong H^p(M;\reals) + H^{p-1}(M;\reals).$$
By definition, we have
$$H^p(\Harm^\dast(M),d) = {\Charm^p\over 
d(\Harm^{p-1})},$$
where $\Charm^p$ denotes the set $C^p\cap \Harm^p$ of $p$-forms
which are both closed and harmonic. 
Recalling that $CcC^p_N$ is the orthogonal complement of the exact
$p$-forms
within the closed ones, we can write
$$\Charm^p = CcC^p_N \oplus \Eharm^p,$$
where $\Eharm^p$ denotes the space of exact harmonic $p$-forms.
We naturally have \\
$d(\Harm^{p-1})\subset 
\Eharm^p$, and thus get a direct-sum
decomposition
$$H^p(\Harm^\dast(M),d) = CcC^p_N + {\Eharm^p
\over d(\Harm^{p-1})}.$$
The first term on the right is isomorphic to $H^p(M;\reals)$. 
The second term on the right measures the extent to which a harmonic
$p$-form can be exact without actually being the exterior derivative
of a harmonic $p-1$-form. This is the term that we claim to be
the echo of $H^{p-1}(M;\reals)$. As suggested by the  
example in section 1, this isomorphism is provided
by the co-differential $\delta$. We demonstrate this in the following
lemma,
which will complete the proof of Theorem 1.

\medskip

\noindent\bf Lemma 3\it.\ \ Under the assumptions of Theorem 1, the
co-differential $\delta\colon \Omega^p\to
\Omega^{p-1}$ induces an isomorphism
$$\bah{\delta}\colon {\Eharm^p
\over d(\Harm^{p-1})}\to H^{p-1}(M;\reals).$$
\rm

\medskip
\noindent That is, the isomorphism $\bah{\delta}$ takes the echo back to
its 
source.

\noindent \it Proof of Lemma 3\rm.\ \ We show that
the linear map 
$$\bah{\delta}\colon {\Eharm^p
\over d(\Harm^{p-1})}\to {C^{p-1}\over
E^{p-1}} \cong H^{p-1}(M;\reals)$$
is well-defined by seeing that
the numerator of the domain of $\bah{\delta}$ maps to the 
numerator of its range, and likewise for the denominators.
First, if $\ph\in \Eharm^p$, then $\ph$ is
an exact, harmonic $p$-form. Being exact, $\ph$ is certainly
closed, hence $\Delta\ph=(\delta d+d\delta)\ph=d\delta\ph=0$. 
Thus $\delta\ph$ is
a closed $p-1$-form. 
Second, if $\ph\in d(\Harm^{p-1})$ is the exterior derivative of
a harmonic $p-1$-form $\beta$, then $\delta\ph = \delta d\beta
=-d\delta\beta$, showing that $\delta\ph$ is an exact $p-1$-form. 
Hence $\bah\delta$ is well-defined.

Next, we show that $\bah\delta$ is one-to-one. To this end, 
suppose that $\ph\in\Eharm^p$ and that
$\delta\ph\in E^{p-1}$. We must show that $\ph\in d(\Harm^{p-1})$.
Since $\ph$ is exact, write $\ph=d\beta$ for $\beta\in \Omega^{p-1}$, and 
note that the Laplacian of $\beta$ is exact, since
$$\Delta\beta=\delta d\beta+d\delta\beta=\delta\ph +d\delta\beta
\in E^{p-1}.$$ Thus $\Delta\beta =d\eta$ for some $p-2$-form $\eta$.
Since
the Laplacian on $p-2$-forms is surjective 
(Lemma 2), we write $\eta=\Delta\sigma$. Then, because
$\Delta\beta=d\eta=d\Delta \sigma=\Delta d\sigma$, we have that
$\beta-d\sigma$ is harmonic.
Finally, writing
$\ph=d(\beta-d\sigma)$ shows that $\ph\in d(\Harm^{p-1})$, as desired.

Finally, to prove that $\bah{\delta}$ is surjective, given 
$\alpha\in C^{p-1}$, we must find an exact harmonic
form $\ph\in \Eharm^p$ such that $\delta\ph - \alpha
\in E^{p-1}$. 
Using the surjectivity of the Laplacian on $p-1$-forms (Lemma 2 again),
we write $\alpha=\Delta\beta$, and then let $\ph=d\beta$. Note that
$\Delta\ph=\Delta d\beta=d\Delta\beta=d\alpha=0$, since $\alpha$ is closed.
Therefore $\ph$ is harmonic, and hence lies in $\Eharm^p$. 
Now, 
$$\delta\ph=\delta d\beta=\Delta\beta-d\delta\beta=\alpha-d\delta\beta,$$
so $\delta\ph - \alpha=-d\delta\beta$,
showing that $\delta\ph-\alpha$ is exact, as desired.

This completes the proof of Lemma 3, and with it, the proof of Theorem 1.

\bigskip

\noindent\bf References\rm

\addtolength{\baselineskip}{-2pt}

\begin{description}
\item[1931] G.\ de\,Rham, {\it Sur l'analysis situs 
des vari\'{e}t\'{e}s \`a  n  
        dimensions}, J.\ Math.\ Pures Appl., \bf 10\rm, 115--200.
 
\item[1933] W.V.D.\ Hodge, {\it A Dirichlet problem for 
harmonic functions with 
     applications
        to analytic varieties}, Proc.\ London Math.\ Soc., 
\bf 36\rm, 257--303.

\item[1940] H.\ Weyl, {\it The method of orthogonal 
projection in potential theory},
        Duke Math.\ J., \bf 7\rm, 411--444.

\item [1941] W.V.D.\ Hodge, {\it The Theory and Applications 
of Harmonic Integrals},
      Cambridge University Press, Cambridge.

\item [1941] A.\ Tucker, {\it  A boundary-value 
theorem for harmonic tensors},        
           Bull.\ AMS, \bf 47\rm, 714.
  
\item [1943] H.\ Weyl, {\it On Hodge's theory 
of harmonic integrals}, Annals of Math., \bf
  44\rm, 1--6.

\item[1946]  P.\ Bidal and G.\ de\,Rham, {\it Les formes 
diff\'erentielles harmoniques},
          Comment.\ Math.\ Helv., \bf 19\rm, 1--49.

\item[1949] K.\ Kodaira, {\it Harmonic fields in 
Riemannian manifolds (Generalized
        potential theory)}, Annals of Math., \bf 50\rm, 587--665.

\item[1952] G.F.D.\ Duff, {\it Differential forms in 
manifolds with boundary}, Annals of
        Math., \bf 56\rm, 115--127.

\item[1952] G.F.D.\ Duff and D.\ Spencer, {\it  Harmonic 
tensors on Riemannian 
        manifolds with boundary}, Annals of Math.\ \bf 
56\rm, 127--156.

\item[1955]  K.O.\ Friedrichs, {\it Differential forms 
on Riemannian manifolds}, Comm.\ Pure 
   Appl.\ Math., \bf 8\rm, 551--590.
 
\item[1955] G.\ de\,Rham, {\it Vari\'et\'es 
Diff\'erentiables}, Hermann, Paris
        (English edition: \it Differentiable Manifolds\rm, 
Grund.\ der Math.\
        Wiss., \bf 266\rm, Springer-Verlag, Berlin 1984).

\item[1956]  P.\ Conner, {\it The Neumann 
problem for differential forms on 
        Riemannian manifolds}, Memoirs of the AMS \bf 20\rm, 
AMS, Providence.
        
\item[1956]  C.B.\ Morrey, {\it  A variational method 
in the theory of harmonic 
      integrals}, II, 
        Amer.\ J.\ Math., \bf 78\rm, 137--170.
           
\item[1957]  N.\ Aronszajn, {\it  A unique continuation 
theorem for solutions of elliptic
          partial differential equations or inequalities of second order},
          J.\ Math.\ Pures et Appl., \bf 36\rm, 235--249.

\item[1962]  N.\ Aronszajn, A.\ Krzywicki and J.\ Szarski, 
{\it  A unique continuation 
      theorem
        for exterior differential forms on 
Riemannian manifolds}, Ar.\ Kat., \bf 4\rm, 
        417--453.

\item[1966]  C.B.\ Morrey, {\it Multiple Integrals 
in the Calculus of Variations}, 
     Grund.\ der 
       Math.\ Wiss., \bf 130\rm, Springer-Verlag, Berlin.
        
\item[1980] J.\ Cheeger, {\it On the Hodge theory of Riemannian
Pseudomanifolds}, Proc.\ Sympos.\ Pure Math., \bf 36\rm, 
(\it Geometry
of the Laplace operator\rm), 91--146.

\item[1988]  J.\ Kazdan, {\it Unique continuation 
in geometry}, Comm.\ Pure Appl.\ Math., \bf
          41\rm, 667--681.

\item[1995] G.\ Schwarz, \it Hodge Decomposition --- A Method
for Solving Boundary Value Problems\rm, Lecture Notes in Math.,
\bf 1607\rm, Springer-Verlag, Berlin.

\item[1999] R.\ Mazzeo and R.\ Melrose, \it Pseudodifferential operators
on manifolds with fibred boundaries\rm, Asian J.\ of Math., \bf 2\rm, 
833--866.
        
\item[2004]  J.\ Parsley, \it The Biot-Savart operator and
electrodynamics on bounded subdomains of the three-sphere\rm,
Ph.D.\ thesis, University of Pennsylvania.

\item[2005] C.\ Epstein, \it Subelliptic Spin$_{C}$ Dirac operators,
I
and II\rm, to appear.

\end{description}

\noindent
Sylvain Cappell: Courant Institute, New York University, 
\it cappell@cims.nyu.edu\rm\\
Dennis DeTurck: University of Pennsylvania, \it deturck@math.upenn.edu\rm\\
Herman Gluck: University of Pennsylvania, \it gluck@math.upenn.edu\rm\\
Edward Y.\ Miller: Polytechnic University, \it emiller@poly.edu\rm

\end{document}